# Newton Like Iterative Method without Derivative for Solving Nonlinear Equations Based on Dynamical Systems

Yonglong Liao [a*] and Limin Cui [b]

ZHIYUAN College, Beijing Institute of Petrochemical Technology, Beijing, China

[a]liaoyonglong@bipt.edu.cn, [b]limincui@bipt.edu.cn

## Abstract

The iterative problem of solving nonlinear equations is studied. A new Newton like iterative method with adjustable parameters is designed based on the dynamic system theory. In order to avoid the derivative function in the iterative scheme, the difference quotient is used instead of the derivative. Different from the existing methods, the difference quotient scheme in this paper has higher accuracy. Thus, the new iterative method is suitable for a wider range of initial values. Finally, several numerical examples are given to verify the practicability and superiority of the method.

## 1. INTRODUCTION

Solving equation is one of the most important problems in engineering calculation. For nonlinear equations $f(x)=0$, because of its importance and the universality of the problems involved, its study has a long history. Classical and commonly used methods include dichotomy, Newton's method and secant method [1,2]. Among them, Newton's method is a general iterative method in many different situations. Newton's method has second-order convergence. Once it converges, it is much faster than the bisection method and secant method. Unfortunately, this method depends on the selection of initial values heavily, and it cannot be guaranteed that it is always convergent for different problems. Therefore, scholars have proposed various improved structures for Newton's method [3-7], aimed to construct an iterative algorithm with a wider scope of application and a larger convergence domain.

The Newton's method requires that $f'(x) \neq 0$, which restricts the use of Newton's method to solve some nonlinear equations. In order to overcome this shortcoming of Newton's method, literature [8] derived a new iterative formula with the help of dynamic system method, maintaining the advantages of Newton's method and removing the requirement $f'(x) \neq 0$. On this basis, reference [9] gives an iterative method without derivative for solving nonlinear equations, which has second-order convergence with adjustable parameters. References [10-11] proposed a predictor corrector iterative scheme for finding the solutions of nonlinear equations by using the improved Newton like method as predictor, and proved at least second-order convergence.

In the further study, although the condition $f'(x) \neq 0$ is unnecessary, the derivative in the iterative scheme makes it difficult to solve and increase the calculation cost greatly in some problems. For this reason, literature [9] uses difference quotient scheme to replace the derivative. In the difference quotient scheme, $f(x_n)$ is used as the step length rather than $x_n - x_{n-1}$, which brings great convenience to the stability proof, but enlarges the error. In order to improve the convergence effect of the iterative scheme further, this paper introduces two different difference quotients, which can effectively reduce the iteration error on the premise of ensuring at least second-order convergence.

The structure of this paper is as follows: the second part is the mathematical description and lemma of the problem; the third part gives the convergence analysis; on this basis, the fourth part shows the effectiveness of the method through several numerical simulations; the last part is a brief conclusion.

## 2. PROPERTIES MATHEMATICAL DESCRIPTION OF PROBLEMS AND NEW ITERATIVE SCHEME

Consider the numerical solution of nonlinear equation
$$f(x) = 0 \tag{1}$$

where, $x \in R$ represents variable, $f(x)$ represents a real valued continuous function within interval $[a,b]$. Let equation (1) has a solution $x^*$ in the interval $[a,b]$, namely, $\exists x^* \in [a,b]$ satisfies $f(x^*) = 0$.

We introduce the dynamic system [8]
$$\begin{cases} \dfrac{dx}{dt} = -\dfrac{f(x)}{\mu f(x) + f'(x)} \\ x(0) = x_0, \quad x_0 \in U(x^*) \end{cases} \tag{2}$$

where $U(x^*)$ represents a neighborhood of $x^*$, and $\mu$ represents the iteration parameter. It is easy to known that the solution of (1) is the equilibrium point of the dynamic system (2), and vice versa.

**Lemma1**[8]: let $x^*$ is the solution of (1), $f'(x^*) \neq 0$, and $f(x)$ has the second derivative in $U(x^*)$, if the initial value $x_0 \in U(x^*)$, then the dynamic system (2) has a unique solution $x = x(t, x_0)$ in $(0, +\infty)$, and moreover
$$\lim_{t \to +\infty} x(t, x_0) = x^* \tag{3}$$

$U(x^*)$ represents the neighborhood of $x^*$.

Lemma 1 makes it possible to solve nonlinear equation (1) by solving dynamic system (2).

Using the Euler method for the differential equation in (2), we have
$$x_{n+1} = x_n - \frac{h f(x_n)}{\mu f(x_n) + f'(x_n)} \tag{4}$$

where $h > 0$ represents the iteration step length. Let $h = 1$, we have
$$x_{n+1} = x_n - \frac{f(x_n)}{\mu f(x_n) + f'(x_n)} \tag{5}$$

In general, the function $f(x)$ is complex, and the calculation of its derivative is much heavy. In order to avoid the demand of derivative values, literature [9] uses difference quotient instead of derivative to obtain the following lemma.

**Lemma2**[9]: The following iterative scheme for nonlinear equation (1) has second-order convergence at least

$$x_{n+1} = x_n - \frac{f^2(x_n)}{\mu f^2(x_n) + f(x_n + f(x_n)) - f(x_n)} \quad (6)$$

In equation (6), $f(x_n)$ is taken as the step length of difference quotient scheme, which facilitates its convergence analysis. However, this difference quotient step length tends to increase the error of derivative approximation. Now, we improve the difference quotient and give the following iterative scheme:

$$x_{n+1} = x_n - \frac{f(x_n)}{\mu f(x_n) + \frac{f(x_n) - f(x_{n-1})}{x_n - x_{n-1}}}$$

After simplification, we have

$$x_{n+1} = x_n - \frac{f(x_n)(x_n - x_{n-1})}{\mu(x_n - x_{n-1})f(x_n) + f(x_n) - f(x_{n-1})}$$

A new iterative scheme for solving nonlinear equation (1) is obtained as follows:

$$x_{n+1} = x_n - \frac{f(x_n)(x_n - x_{n-1})}{\mu(x_n - x_{n-1})f(x_n) + f(x_n) - f(x_{n-1})} \quad (7)$$

**Remark1**: The difference quotient step length in this paper is $x_n - x_{n-1}$, which can reduce the error of derivative approximation effectively. It will be verified in MATLAB experiment in section 4.

## 3. CONVERGENCE ANALYSIS

In this section, the following second-order convergence and convergence rate estimates are given.

**Theorem1**: let $x^*$ is the solution of (1), $f'(x^*) \neq 0$, and $f(x)$ has the second derivative in $U(x^*)$, if $x_0 \in U(x^*)$, then the iterative scheme (7) has second-order convergence at least, and

$$\lim_{n \to \infty} \frac{x_{n+1} - x^*}{(x_n - x^*)^2} = \mu + \frac{f''(x^*)}{f'(x^*)} \quad (8)$$

Proof: We denote $\overset{\circ}{U}(x^*)$ as the noncentral neighborhood of $x^*$. Let $x_n \in \overset{\circ}{U}(x^*)$ and $e_n = x_n - x^*$, then

$$\begin{aligned}
e_{n+1} &= x_{n+1} - x^* \\
&= x_n - \frac{f(x_n)(x_n - x_{n-1})}{\mu(x_n - x_{n-1})f(x_n) + f(x_n) - f(x_{n-1})} - x^* \\
&= e_n - \frac{f(x_n)(x_n - x_{n-1})}{\mu(x_n - x_{n-1})f(x_n) + f(x_n) - f(x_{n-1})}
\end{aligned} \quad (9)$$

Considering the dynamic (2), the backward differential is used at $x_n$, we have

$$\frac{x_n - x_{n-1}}{h} = -\frac{f(x_n)}{\mu f(x_n) + f'(x_n)}$$

By taking $h=1$, it has the form

$$x_{n-1} = x_n + \frac{f(x_n)}{\mu f(x_n) + f'(x_n)} \quad (10)$$

Substitute (10) into (9) to get

$$e_{n+1} = e_n + \frac{f^2(x_n)}{f(x_n)f'(x_n) - [\mu f(x_n) + f'(x_n)]f(x_n + \frac{f(x_n)}{\mu f(x_n) + f'(x_n)})} \quad (11)$$

Taking the Taylor expansion of $f(x_n)$ and $f'(x_n)$, we have

$$f(x_n) = f(x^* + e_n) = f'(x^*)e_n + \frac{1}{2}f''(x^*)e_n^2 + o(e_n^2) \quad (12)$$

$$f'(x_n) = f'(x^*) + f''(x^*)e_n + \frac{1}{2}f'''(x^*)e_n^2 + o(e_n^2) \quad (13)$$

Combining (12) and (13), we get

$$\frac{f(x_n)}{\mu f(x_n) + f'(x_n)}$$

$$= \frac{f'(x^*)e_n + \frac{1}{2}f''(x^*)e_n^2 + o(e_n^2)}{f'(x^*) + \mu f'(x^*)e_n + f''(x^*)e_n + \frac{1}{2}\mu f''(x^*)e_n^2 + \frac{1}{2}f'''(x^*)e_n^2 + o(e_n^2)}$$

$$= o(e_n) \quad (14)$$

Taking the Taylor expansion of $f(x_n + \frac{f(x_n)}{\mu f(x_n) + f'(x_n)})$, and combining formula (14), we have

$$f(x_n + \frac{f(x_n)}{\mu f(x_n) + f'(x_n)})$$

$$= f'(x^*)\left[e_n + \frac{f(x_n)}{\mu f(x_n) + f'(x_n)}\right] + \frac{1}{2}f''(x^*)\left[e_n + \frac{f(x_n)}{\mu f(x_n) + f'(x_n)}\right]^2$$

$$+ o\left[e_n + \frac{f(x_n)}{\mu f(x_n) + f'(x_n)}\right]^2$$

$$= f'(x^*)e_n + \frac{1}{2}f''(x^*)e_n^2 + \frac{f'(x^*)f(x_n)}{\mu f(x_n) + f'(x_n)} + \frac{f''(x^*)f(x_n)e_n}{\mu f(x_n) + f'(x_n)}$$

$$+ \frac{f''(x^*)f^2(x_n)}{2[\mu f(x_n) + f'(x_n)]^2} + o(e_n^2)$$

Then, we have

$$[\mu f(x_n) + f'(x_n)]f(x_n + \frac{f(x_n)}{\mu f(x_n) + f'(x_n)})$$

$$= \left[f'(x^*)e_n + \frac{1}{2}f''(x^*)e_n^2\right][\mu f(x_n) + f'(x_n)] + f'(x^*)f(x_n) + f''(x^*)f(x_n)e_n$$

$$+\frac{f''(x^*)f^2(x_n)}{2[\mu f(x_n)+f'(x_n)]}+o(e_n^2)$$

$$=2f'^2(x^*)e_n+3f'(x^*)f''(x^*)e_n^2+\mu f'^2(x^*)e_n^2$$

$$+\frac{f'^2(x^*)f''(x^*)e_n^2+f'(x^*)f''^2(x^*)e_n^3+o(e_n^3)}{2f'(x^*)+2\mu f'(x^*)e_n+2f''(x^*)e_n+\mu f''(x^*)e_n^2+f'''(x^*)e_n^2+o(e_n^2)}+o(e_n^2)$$

Therefore, the denominator of the second part of equation (11) can be written as

$$\phi = f(x_n)f'(x_n)-[\mu f(x_n)+f'(x_n)]f(x_n+\frac{f(x_n)}{\mu f(x_n)+f'(x_n)})$$

$$=-f'^2(x^*)e_n-\frac{3}{2}f'(x^*)f''(x^*)e_n^2-\mu f'^2(x^*)e_n^2$$

$$-\frac{f'^2(x^*)f''(x^*)e_n^2+f'(x^*)f''^2(x^*)e_n^3+o(e_n^3)}{2f'(x^*)+2\mu f'(x^*)e_n+2f''(x^*)e_n+\mu f''(x^*)e_n^2+f'''(x^*)e_n^2+o(e_n^2)}+o(e_n^2) \quad (15)$$

Thus, it can be calculated that

$$\varphi = \phi e_n + f^2(x_n)$$

$$=-f'^2(x^*)e_n^2-\frac{3}{2}f'(x^*)f''(x^*)e_n^3-\mu f'^2(x^*)e_n^3$$

$$-\frac{f'^2(x^*)f''(x^*)e_n^3+f'(x^*)f''^2(x^*)e_n^4+o(e_n^4)}{2f'(x^*)+2\mu f'(x^*)e_n+2f''(x^*)e_n+\mu f''(x^*)e_n^2+f'''(x^*)e_n^2+o(e_n^2)}+o(e_n^3)$$

$$+f'^2(x^*)e_n^2+f'(x^*)f''(x^*)e_n^3+o(e_n^3)$$

$$=-\frac{1}{2}f'(x^*)f''(x^*)e_n^3-\mu f'^2(x^*)e_n^3$$

$$-\frac{f'^2(x^*)f''(x^*)e_n^3+f'(x^*)f''^2(x^*)e_n^4+o(e_n^4)}{2f'(x^*)+2\mu f'(x^*)e_n+2f''(x^*)e_n+\mu f''(x^*)e_n^2+f'''(x^*)e_n^2+o(e_n^2)}+o(e_n^3) \quad (16)$$

Therefore, formula (11) can be written as

$$e_{n+1}=\frac{\varphi}{\phi} \quad (17)$$

where $\phi$ and $\varphi$ are given by (15) and (16).

According to formula (17), we have

$$\lim_{n\to\infty}\frac{e_{n+1}}{e_n^2}=\frac{-\frac{1}{2}f'(x^*)f''(x^*)-\mu f'^2(x^*)-\frac{f'^2(x^*)f''(x^*)}{2f'(x^*)}}{-f'^2(x^*)}=\frac{\mu f'(x^*)+f''(x^*)}{f'(x^*)} \quad (18)$$

Namely, the iterative scheme (7) has second-order convergence at least, and the result (8) is obtained.

Theorem 1 shows that the iterative scheme (7) is at least second-order convergent. When appropriate parameters are obtained, such as $\mu=-\frac{f''(x^*)}{f'(x^*)}$, the convergence rate can be higher than the second order.

## 4. NUMERICAL SIMULATIONS

In this section we introduce the following three examples to verify the conclusions of this paper.

**Example1** (Logarithmic equation) $f(x) = \ln x$, the range of variable $x$ is $[a,b] = [0.5, 5]$, the initial value is $x_0 = 5$.

**Example2** (Exponential equation) $f(x) = (x-1)e^{-x}$, the range of variable $x$ is $[a,b] = [-1, 50]$, the initial value is $x_0 = 50$.

**Example3** (Trigonometric function equation) $f(x) = 2\sin x - 1$, the range of variable $x$ is $[a,b] = [0, \frac{11\pi}{24}]$, the initial value is $x_0 = \frac{11\pi}{24}$.

The calculation results of Example 1-3 are shown in Table 1. The iterative scheme of Newton's method is

$$x_{n+1} = x_n - \frac{f(x_n)}{f'(x_n)}$$

The parameters in iteration scheme (6) are taken as $\mu = 1$, $\mu = 1 + \frac{1}{e}$, $\mu = \frac{1}{2} + \frac{\sqrt{3}}{6}$ in examples 1-3 respectively, which are consistent with the values in reference [9], and are the best parameter values of the iterative scheme. The parameters in iteration scheme (7) are taken as $\mu = 0.135$, $\mu = 1.18$, $\mu = 2.65$ in examples 1-3 respectively. The maximum iterations are set to 500, and the precision is taken as $e = 10^{-5}$.

**Table 1.** Convergence and numerical solutions of iterative schemes

| example | initial value | exact solution | Newton's method | scheme (6) | scheme (7) | |
|---|---|---|---|---|---|---|
| | | | iteration times | iteration times | iteration times | $x_n$ |
| 1 | 5 | 1 | divergence | $n = 8$ | $n = 6$ | 1.000000 |
| 2 | 50 | 1 | divergence | divergence | $n = 40$ | 1.000000 |
| 3 | $\frac{11\pi}{24}$ | $\frac{\pi}{6}$ | divergence | divergence | $n = 6$ | 0.523599 |

It can be seen from Table 1 that the results are divergent by using the classical Newton's method. Only the result of Example 1 is convergent by the iterative scheme (6), and the other two examples are divergent. All of the three examples are convergent by using the iterative scheme (7) given in this paper, and the convergence rate in Example 1 is significantly faster than that of the iterative scheme (6).

## 5. CONCLUSION

This paper presents a new iterative method for solving nonlinear equations. The method has adjustable parameters and does not require derivative function. Compared with the methods in the existing literature, a larger initial value selection range is allowed. The method is proved to be at least second-order convergent by theoretical derivation. Finally, several numerical examples are given to demonstrate the superiority of this method.

In the design process, we find that the iteration step length $h$ has some influence on the convergence rate. If the value of $h$ is smaller, the iteration result of each step is more accurate, but the iteration time is prolonged. When the value of $h$ is larger, the iteration time will be shortened, but the iteration results of each step are rough. How to select an appropriate iteration step length $h$ is the next problem we will study.

## ACKNOWLEDGEMENTS

This paper was supported by the general project of science and technology plan of Beijing Municipal Commission of Education (No. KM202110017002); curriculum development and teaching reform project of Beijing Institute of Petrochemical Technology (No. YJ22-202).

## REFERENCES

[1] R.L. Burden, J.D. Faires and A.M. Burden: *Numerical Analysis*, 10th edition (Cengage Learning: Boston, MA, USA, 2016).

[2] D. Kincaid, W. Cheney: *Numerical Analysis: Mathematics of Scientific Computing*, 3rd edition (American Mathematical Society, Providence, Rhode Island, USA, 2002).

[3] S. Weerakoon, T.G.I. Fernando: A variant of Newton's method with accelerated third-order convergence. Applied Mathematics Letters, Vol. 13 (2000) No. 8, p. 87-93.

[4] A.Y. Özban: Some new variants of Newton's method, Applied Mathematics Letters, Vol. 17 (2004) No. 6, p. 667-682.

[5] M. Frontini, E. Sormani: Some variant of Newton's method with third-order convergence, Applied Mathematics and Computation, Vol. 140 (2003) No. 2-3, p. 419-426.

[6] H.H.H. Homeier: On Newton-type methods with cubic convergence, Journal of Computational and Applied Mathematics, Vol. 176 (2005) No. 2, p. 425-432.

[7] T. Lukic, N.M. Ralevic: Geometric mean Newton's method for simple and multiple roots, Applied Mathematics Letters, Vol. 21(2008) No. 1, p. 30-36.

[8] X. Wu: A significant improvement on newton's iterative method, Applied Mathematics and Mechanics, Vol. 20(1999) No. 8, p. 863-866. (In Chinese)

[9] Q. Zheng: Parametric iterative methods of quadratic convergence without the derivative, Mathematica Numerica Sinica, Vol. 25(2003) No. 1, p.107-112. (In Chinese)

[10] H. Cai and L. Qian: Predict-correct iterative methods for solving foots of nonlinear equation, Journal of Yili Normal University（Natural Science Edition）, Vol. 1(2010) No.1, p. 17-19. (In Chinese)

[11] K. Rahman: One class of third-order iteration methods for solving non-linear equations, Journal of Jiangxi Normal University ( Natural Science), Vol. 44(2020) No. 2, p. 206-208. (In Chinese)